\documentclass[12pt,twoside,onecolumn,english]{svjour3}
\usepackage[T1]{fontenc}
\usepackage[latin9]{inputenc}
\usepackage{amsmath}
\usepackage{amssymb}
\usepackage{graphicx}

\makeatletter

\providecommand{\tabularnewline}{\\}

\smartqed  

\authorrunning{Y. Censor, R. Davidi, G.T. Herman, R.W. Schulte and L. Tetruashvili}
\titlerunning{Projected subgradient minimization versus superiorization}

\@ifundefined{showcaptionsetup}{}{%
 \PassOptionsToPackage{caption=false}{subfig}}
\usepackage{subfig}
\makeatother

\usepackage{babel}
\begin{document}

\title{Projected Subgradient Minimization versus Superiorization }

\date{February 5, 2013. Revised: June 11, 2013, and August 14, 2013. \\
Communicated by Masao Fukushima}

\author{Yair Censor$^{*1}$, Ran Davidi$^{2}$, Gabor T. Herman$^{3}$, Reinhard
W. Schulte$^{4}$ and Luba Tetruashvili$^{1}$}

\institute{$^{1}$Department of Mathematics, University of Haifa, Mt. Carmel,
3190501 Haifa, Israel ($^{*}$corresponding author: Y. Censor, e-mail:
yair@math.haifa.ac.il) \\
$^{2}$Department of Radiation Oncology, Stanford University, Stanford,
CA 94305, USA\\
$^{3}$Department of Computer Science, The Graduate Center, City University
of New York, New York, NY 10016, USA \\
$^{4}$Department of Radiation Medicine, Loma Linda University Medical
Center, Loma Linda, CA 92354, USA }
\maketitle
\begin{abstract}
The projected subgradient method for constrained minimization repeatedly
interlaces subgradient steps for the objective function with projections
onto the feasible region, which is the intersection of closed and
convex constraints sets, to regain feasibility. The latter poses a
computational difficulty and, therefore, the projected subgradient
method is applicable only when the feasible region is ``simple to
project onto\textquotedblright{}. In contrast to this, in the superiorization
methodology a feasibility-seeking algorithm leads the overall process
and objective function steps are interlaced into it. This makes a
difference because the feasibility-seeking algorithm employs projections
onto the individual constraints sets and not onto the entire feasible
region.

We present the two approaches side-by-side and demonstrate their performance
on a problem of computerized tomography image reconstruction, posed
as a constrained minimization problem aiming at finding a constraint-compatible
solution that has a reduced value of the total variation of the reconstructed
image.


\end{abstract}

\section{Introduction\label{sect:introduction-1}}

Our aim in this paper is to expose the recently-developed superiorization
methodology and its ideas to the optimization community by ``confronting''
it with the projected subgradient method. We juxtapose the\textit{
projected subgradient} \textit{method }(PSM) with the \textit{superiorization
methodology} (SM) and demonstrate their performance on a large-size
real-world application that is modeled, and needs to be solved, as
a constrained minimization problem. The PSM for constrained minimization
has been extensively investigated, see, e.g., \cite[Subsection 7.1.2]{AR06},
\cite[Subsection 3.2.3]{nestrov-book}. Its roots are in the work
of Shor \cite{shor-book} for the unconstrained case and in\textbf{
}the work of Polyak \cite{polyak67,polyak69} for the constrained
case. More recent work can be found in, e.g., \cite{bt03}. The superiorization
methodology was first proposed in \cite{bdhk07}, although without
using the term superiorization. In that work, perturbation resilience
(without using this term) was proved for the general class of \textit{string-averaging}
\textit{projection }(SAP)\textit{ }methods, see \cite{ceh01,CS08,CS09,ct03,pen09},
that use orthogonal projections and relate to consistent constraints.
Subsequent investigations and developments of the SM were done in
\cite{cdh10,dhc09,hd08,ndh12,pscr10}. More information on superiorization-related
work is given in Section \ref{sec:related-work}. 

It is not claimed that the PSM is the best optimization method for
solving constrained minimization problems and there are many different
alternative methods with which SM could be compared. So, why did we
chose to confront the PSM with our SM? In a nutshell, our answer is
that both methods interlace steps related to the objective function
with steps oriented toward feasibility, but they differ in how they
restore or preserve feasibility.A major difficulty with the PSM is
the need to perform, within each iterative step, an orthogonal projection
onto the feasible set of the constrained minimization problem. If
the feasible set is not ``simple to project onto\textquotedblright{}\ then
the projection requires an independent inner-loop calculation to minimize
the distance from a point to the feasible set, which can be costly
and hamper the overall effectiveness of the PSM.

In the SM, we replace the notion of a fixed feasible set by that of
a nonnegative real-valued proximity function. This function serves
as an indicator of how incompatible a vector is with the constraints.
In such a formulation, the merit of an actual output vector of any
algorithm is indicated by the smallness of the two numbers, i.e.,
the values of the proximity function and the objective function.The
underlying idea of SM is that many iterative algorithms that produce
outputs for which the proximity function is small are \textit{strongly
perturbation resilient} in the sense that, even if certain kinds of
changes are made at the end of each iterative step, the algorithm
still produces an output for which the proximity function is not larger.
This property is exploited by using permitted changes to steer the
algorithm to an output that has not only a small proximity function
value, but has also a small objective function value.

The PSM requires that feasibility is regained after each subgradient
step by performing a projection onto the entire feasible set whereas
in the SM the feasibility-seeking projection method proceeds by projecting
(in a well-defined algorithmically-structured regime dictated by the
specific projection method) onto the individual sets, whose intersection
is the entire feasible set, and not onto the whole feasible set itself.
This has a potentially great computational advantage.

We elaborate on the motivation for this work in Section \ref{sect:Motivation}.
In Section \ref{sec:related-work} we discuss some superiorization-related
work, in Section \ref{sect:SM} the SM is presented, and in Section
\ref{sect:experimental} we demonstrate the approaches of the SM and
the PSM on a realistically-large-size problem with data that arise
from the significant problem of x-ray computed tomography (CT) with
total variation (TV) minimization, followed by some conclusions in
Section \ref{sect:conclusions}.

\section{Motivation and Basic Notions\label{sect:Motivation}}

Throughout this paper, we assume that $\Omega$ is a nonempty subset
of the $J$-dimensional Euclidean space $\mathbb{R}^{J}$. We consider
constrained minimization problems of the form 
\begin{equation}
\mathrm{minimize}\left\{ \phi(x)\mid x\in C\right\} ,\label{prob:cons-min}
\end{equation}
where $\phi:\mathbb{R}^{J}\rightarrow\mathbb{R}$ is an objective
function and $C\subseteq\Omega$ is a given feasible set. 

Since we juxtapose the\textit{ projected subgradient} \textit{method
}(PSM) with the \textit{superiorization methodology} (SM) and demonstrate
their performance on a large-size real-world application that is modeled,
and needs to be solved, as a constrained minimization problem, we
now outline these two methods and explain our choice in detail.

In order to apply the PSM to solving (\ref{prob:cons-min}) we need
to assume that $C$ is a nonempty closed convex set and that $\phi$
is a convex function. The PSM generates a sequence of iterates $\left\{ x^{k}\right\} _{k=0}^{\infty}$
according to the recursion formula
\begin{equation}
x^{k+1}=P_{C}\left(x^{k}-t_{k}\phi^{\prime}\left(x^{k}\right)\right),\label{eq:sgp}
\end{equation}
where $t_{k}>0$ is a step-size, $\phi^{\prime}\left(x^{k}\right)\in\partial\phi\left(x^{k}\right)$
is a subgradient of $\phi$ at $x^{k},$ and $P_{C}$ stands for the
orthogonal (least Euclidean norm) projection onto the set $C.$ 

A major difficulty with (\ref{eq:sgp}) is the need to perform, within
each iterative step, the orthogonal projection. If the feasible set
$C$ is not ``simple to project onto\textquotedblright{}\ then the
projection requires an independent inner-loop calculation to minimize
the distance from the point $x^{k}-t_{k}\phi^{\prime}\left(x^{k}\right)$
to the set $C$, which can be costly and hamper the overall effectiveness
of an algorithm that uses (\ref{eq:sgp}). Also, if the inner loop
converges to the projection onto $C$ only in the limit, then, in
practical implementations, it will have to be stopped after a finite
number of steps, and so $x{}^{k+1}$ will be only an approximation
to the projection onto $C$ and it could even happen that it is not
in $C$.

Even if we set aside our worries about projecting onto $C$ in (\ref{eq:sgp}),
there are still two concerns when applying the PSM to real-world problems.
One is that the iterative process usually converges to the desired
solution only in the limit. In practice, some stopping rule is applied
to terminate the process and the output at that time may not even
be in $C$ and, even if it is in $C$, it is most unlikely to be the
minimizer of $\phi$ over $C$. The second problem in real-world applications
comes from the fact that the constraints, derived from the real-world
problem, may not be consistent (e.g., because they come from noisy
measurements) and so $C$ is empty.

Similar criticism applies actually to many constrained-minimization-seeking
algorithms for which asymptotic convergence results are available.
In the SM, both of these objections can be handled by replacing the
notion of a fixed feasible set $C$ by that of a nonnegative real-valued
proximity function $Prox_{C}:\Omega\rightarrow\mathbb{R}_{+}$. This
function serves as an indicator of how incompatible a vector $x$
is with the constraints. In such a formulation, the merit of the actual
output $x$ of any algorithm is indicated by the smallness of the
two numbers $Prox_{C}(x)$ and $\phi(x)$. For the formulation of
(\ref{prob:cons-min}), we would define $Prox_{C}$ so that its range
is the ray of nonnegative real numbers with $Prox_{C}(x)=0$ if, and
only if, $x\in C$ and then the constrained minimization problem (\ref{prob:cons-min})
is precisely that of finding an $x$ that is a minimizer of $\phi(x)$
over $\left\{ x\mid Prox_{C}(x)=0\right\} $. The above discussion
allows us to do away with the nonemptiness assumption and also to
compare the merits of actual outputs of algorithms that only approximate
the aim of the constrained minimization problem.

The recently invented SM incorporates the ideas of the previous paragraph
in its very foundation and formulates the problem with the function
$Prox_{C}$ instead of the set $C$. The underlying idea of SM is
that many iterative algorithms that produce outputs $x$ for which
$Prox_{C}(x)$ is small are \textit{strongly perturbation resilient}
in the sense that, even if certain kinds of changes are made at the
end of each iterative step, the algorithm still produces an output
$x^{\prime}$ for which $Prox_{C}(x^{\prime})$ is not larger. This
property is exploited by using permitted changes to steer the algorithm
to an output that has not only a small $Prox_{C}$ value, but has
also a small $\phi$ value. The algorithm that incorporates such a
steering process is referred to as the \textit{superiorized version}
of the original iterative algorithm. The main practical contribution
of SM is the automatic creation of the superiorized version, according
to a given objective function $\phi$, of just about any iterative
algorithm that aims at producing an $x$ for which $Prox_{C}(x)$
is small.

Nevertheless, in order to carry out our comparative study, we restrict
our attention here to a subset of all possible problems to which not
only the SM but also the PSM is applicable. We assume that we are
given a family of constraints $\left\{ C_{\ell}\right\} _{\ell=1}^{L}$,
where each set $C_{\ell}$ is a nonempty closed convex subset of $\mathbb{R}^{J}$
such that
\begin{equation}
C=\bigcap_{\ell=1}^{L}C_{\ell}\label{eq:inersection}
\end{equation}
is a nonempty subset of $\Omega$ and that it is the feasible set
$C$ of (\ref{prob:cons-min}). Under these assumptions, we illustrate
the application of the SM by the superiorization of feasibility-seeking
\textit{projection method}s, see, e.g., \cite{ac89,bb96,BC11,cccdh10,CZ97}
and the recent monograph \cite{Ceg-book}. Such methods use projections
onto the individual sets $C_{\ell}$ in order to generate a sequence
$\left\{ x^{k}\right\} _{k=0}^{\infty}$ that converges to a point
$x^{\ast}\in C$. Therefore, contrary to the PSM, one does not need
to assume that $C$ is a ``simple to project onto\textquotedblright{}\ set,
but rather that the individual sets $C_{\ell}$ have this property.
The latter is indeed often the case, such as, for example, when the
sets $C_{\ell}$ are hyperplanes or half-spaces onto which we can
project easily, but their intersection is not ``simple to project
onto\textquotedblright{}.

The SM is accurately presented in Section \ref{sect:SM} below. However,
the discussion above is sufficient to explain why we chose the PSM
and the SM for our comparative study. Namely, both methods interlace
objective-function-reduction steps with steps oriented toward feasibility.
But exactly here lies a big difference between the two approaches.
The PSM requires that feasibility is regained after subgradient nonascent
steps by performing a projection onto $C$, whereas in the SM the
feasibility-seeking projection method proceeds by projecting (in a
well-defined algorithmically-structured regime dictated by the specific
projection method) onto the individual sets $C_{\ell}$ and not onto
the whole feasible set $C.$ This has a potentially great computational
advantage.

\section{Superiorization-Related Previous Work\label{sec:related-work}}

The superiorization methodology was first proposed in \cite{bdhk07},
although without using the term superiorization. In that work, perturbation
resilience (without using this term) was proved for the general class
of \textit{string-averaging} \textit{projection }(SAP)\textit{ }methods,
see \cite{ceh01,CS08,CS09,ct03,pen09}, that use orthogonal projections
and relate to consistent constraints. Subsequent investigations and
developments were done in \cite{cdh10,dhc09,hd08,ndh12,pscr10}. In
\cite{cdh10}, the methodology was formulated over general \textit{problem
structures} which enabled rigorous analysis and revealed that the
approach is not limited to feasibility and optimization. In \cite{dhc09},
perturbation resilience was analyzed for the class of \textit{block-iterative
projection }(BIP)\textit{ }methods, see \cite{ac89,bb96,BC11,cccdh10,CZ97},\textit{
}and applied in this manner. In \cite{hd08}, the advantages of superiorization
for image reconstruction from a small number of projections was studied,
and in \cite{ndh12} two acceleration schemes based on (symmetric
and nonsymmetric) BIP methods were proposed and experimented with.
In \cite{pscr10}, total variation superiorization schemes in proton
computed tomography (pCT) image reconstruction were investigated.

In \cite{hgdc12}, we introduced the notion of $\varepsilon$-compatibility
into the superiorization approach in order to handle inconsistent
constraints. This enabled us to close the logical discrepancy between
the assumption of consistency of constraints and the actual experimental
work done previously. We also introduced there the new notion of strong
perturbation resilience which generalizes the previously used notion
of perturbation resilience. Algorithmically, the new superiorized
algorithm introduced there (and used here) is different from all previous
ones in that it uses the notion of nonascending direction and in that
it allows several perturbation steps for each feasibility-seeking
step, an aspect that has practical advantages.

In \cite{jcj12}, superiorization was applied to the \textit{expectation
maximization} (EM) algorithm instead of the feasibility-seeking projection
methods that were used in superiorization previously. The approach
was implemented there to solve an inverse problem of \textit{bioluminescence
tomography} (BLT) image reconstruction. Such EM superiorization was
investigated further and applied to a problem of \textit{Single Photon
Emission Computed Tomography} (SPECT) in \cite{lz12}. Most recently,
in \cite{bk12}, the SM was further investigated numerically, along
with many projection methods for the feasibility problem and for the
best approximation problem.

Our superiorization methodology should be distinguished from the works
of Helou Neto and De Pierro \cite{hdp-siam09,hdp11-1}, of Nedi\'{c}
\cite{Nedic2011}, Ram, Nedi\'{c} and Veeravalli \cite{Nedic2009},
and of Nurminski \cite{nur08b,nur08a,nur10,nur11}. The lack of cross-referencing
between some of these papers shows that, in spite of the similarities
between their approaches, their results were apparently reached independently. 

There are various differences among the works mentioned in the previous
paragraph, differences in overall setup of the problems, differences
in the assumptions used for the various convergence results, etc.
This is not the place for a full review of all these differences.
But we wish to clarify the fundamental difference between them and
the SM. The point is that when two activities are interlaced, here,
feasibility steps and objective function reduction steps, then once
the process is running all such methods look alike. From looking at
the iterative formulas, one cannot tell if (a) ``feasibility steps
are interlaced into an iterative gradient scheme for objective function
minimization'' or if (b) ``objective function reduction steps are
interlaced into an iterative projections scheme for feasibility-seeking''.
The common thread of all works mentioned in the previous paragraph
is that they fall into the category (a), while the SM is of the kind
(b). In all methods of category (a) the condition that is needed to
guarantee convergence to a constrained minimum point is that the diminishing
step-sizes $\alpha_{k}\rightarrow0$ as $k\rightarrow\infty$ must
be such that $\sum_{k=0}^{\infty}\alpha_{k}=+\infty.$ In contrast,
since the feasibility-seeking projection method is the ``leader''
of the overall process in the SM, we must have that the perturbations
(that do the objective function reduction) will use diminishing step-sizes
$\beta_{k}\rightarrow0$ as $k\rightarrow\infty$ but such that $\sum_{k=0}^{\infty}\beta_{k}<\infty.$
The latter condition guarantees the perturbation resilience of the
original feasibility-seeking projection method so that, regardless
of the interlaced objective function reduction steps, the overall
process converges to a feasible, or $\varepsilon$-compatible, point
of the constraints. 

Yet another fundamental difference between the superiorization methodology
and the algorithms of category (a) mentioned above is that those algorithms
perform the interlaced objective function descent and feasibility
steps alternatingly according to a rigid predetermined scheme, whereas
in the superiorization methodology the activation of these steps and
the decisions whether to keep an iterate or discard it are done inside
the superiorized algorithm in a controlled and automatically-supervised
manner. Thus, the superiorization methodology has the following features
not present in the algorithms of category (a) mentioned above: (i)
it conducts iterations of a feasibility-seeking projection method
which is strongly perturbation resilient (as defined below), (ii)
it interlaces objective function nonascent steps into the process
in a controlled and automatically-supervised manner, (iii) it is not
known to guarantee convergence to a solution of the constrained minimization
problem, and it might (we do not know if this is so or not) instead
only be shown to lead to a feasible point whose objective function
value is less than that of a feasible point that would have been reached
by the same feasibility-seeking projection method without the perturbations
exercised by the superiorized algorithm.

The \textit{adaptive steepest descent projections onto convex sets}
(ASD-POCS) algorithm described in \cite{sidk08} has some similarities
to the SM. However, it is not as general as the SM; see \cite{hgdc12}
for a comparison.

\section{The Superiorization Methodology\label{sect:SM}}

In this section we present a restricted version of the SM of \cite{hgdc12}
adapted to our problem (\ref{prob:cons-min}). As discussed in Section
\ref{sect:Motivation}, we associate with the feasible set $C$ in
(\ref{prob:cons-min}) a proximity function ${Prox}_{C}:\Omega\rightarrow\mathbb{R}_{+}$
that is an indicator of how incompatible an $x\in\Omega$ is with
the constraints. For any given $\varepsilon>0$, a point $x\in\Omega$
for which ${Prox}_{C}(x)\leq\varepsilon$ is called an $\varepsilon$\textit{-compatible
solution} for $C$. We further assume that we have, for the $C$ in
(\ref{prob:cons-min}), a feasibility-seeking \textit{algorithmic
operator} $\boldsymbol{A}_{C}:\mathbb{R}^{J}\rightarrow\Omega$, with
which we define the following basic algorithm.

\medskip{}

\textbf{The Basic Algorithm}\\
(B1) \textbf{Initialization}: Choose an arbitrary $x^{0}\in\Omega$,\\
(B2) \textbf{Iterative Step}: Given the current iterate $x^{k}$,
calculate the next iterate $x^{k+1}$ by
\begin{equation}
x^{k+1}=\boldsymbol{A}_{C}\left(x^{k}\right).\label{-The-Basic}
\end{equation}
\medskip{}

The following definition helps to evaluate the output of the Basic
Algorithm upon termination by a stopping rule.\\
\\
\textbf{Definition 4.1} \textbf{The }$\varepsilon$\textbf{-output
of a sequence }

Given $C\subseteq\mathbb{R}^{J}$, a proximity function ${Prox}_{C}:\Omega\rightarrow\mathbb{R}_{+}$,
a sequence $\left\{ x^{k}\right\} _{k=0}^{\infty}\subset\Omega$ and
an $\varepsilon>0,$ then an element $x^{K}$ of the sequence which
has the properties: (i) ${Prox}_{C}\left(x^{K}\right)\leq\varepsilon,$
and (ii) ${Prox}_{C}\left(x^{k}\right)>\varepsilon$ for all $0\leq k<K,$
is called an $\varepsilon$\texttt{-output of the sequence }$\left\{ x^{k}\right\} _{k=0}^{\infty}$\texttt{
with respect to the pair }$(C,$\texttt{ }${Prox}_{C})$. We denote
it by $O\left(C,\varepsilon,\left\{ x^{k}\right\} _{k=0}^{\infty}\right)=x^{K}.$
\\

Clearly, an $\varepsilon$-output $O\left(C,\varepsilon,\left\{ x^{k}\right\} _{k=0}^{\infty}\right)$
of a sequence $\left\{ x^{k}\right\} _{k=0}^{\infty}$ might or might
not exist, but if it does, then it is unique. If $\left\{ x^{k}\right\} _{k=0}^{\infty}$
is produced by an algorithm intended for the feasible set $C,$ such
as the Basic Algorithm, without a termination criterion, then $O\left(C,\varepsilon,\left\{ x^{k}\right\} _{k=0}^{\infty}\right)$
is the \textit{output} produced by that algorithm when it includes
the termination rule to stop when an $\varepsilon$-compatible solution
for $C$ is reached.\\
\\
\textbf{Definition 4.2 Strong perturbation resilience}

Assume that we are given a $C\subseteq\Omega$, a proximity function
${Prox}_{C}$, an algorithmic operator $\boldsymbol{A}_{C}$ and an
$x^{0}\in\Omega$. We use $\left\{ x^{k}\right\} _{k=0}^{\infty}$
to denote the sequence generated by the Basic Algorithm when it is
initialized by $x^{0}$. The Basic Algorithm is said to be\texttt{
strongly perturbation resilient} iff the following hold:

(i) there exist an $\varepsilon>0$ such that the $\varepsilon$-output
$O\left(C,\varepsilon,\left\{ x^{k}\right\} _{k=0}^{\infty}\right)$
exists for every $x^{0}\in\Omega$;

(ii) for every $\varepsilon>0,$ for which the $\varepsilon$-output
$O\left(C,\varepsilon,\left\{ x^{k}\right\} _{k=0}^{\infty}\right)$
exists for every $x^{0}\in\Omega$, we have also that the $\varepsilon^{\prime}$-output
$O\left(C,\varepsilon^{\prime},\left\{ y^{k}\right\} _{k=0}^{\infty}\right)$
exists for every $\varepsilon^{\prime}>\varepsilon$ and for every
sequence $\left\{ y^{k}\right\} _{k=0}^{\infty}$ generated by
\begin{equation}
y^{k+1}=\boldsymbol{A}_{C}\left(y^{k}+\beta_{k}v^{k}\right),\text{ for all }k\geq0,\label{eq:perturb}
\end{equation}
where the vector sequence $\left\{ v^{k}\right\} _{k=0}^{\infty}$
is bounded and the scalars $\left\{ \beta_{k}\right\} _{k=0}^{\infty}$
are such that $\beta_{k}\geq0$, for all $k\geq0,$ and $\sum_{k=0}^{\infty}\beta_{k}<\infty$.
\\
\\
\textbf{Definition 4.3} \textbf{Bounded convergence}Assume that we
are given a $C\subseteq\mathbb{R}^{J}$, a proximity function ${Prox}_{C}$
and an algorithmic operator $\boldsymbol{A}_{C}:\mathbb{R}^{J}\rightarrow\Omega$.
Then the Basic Algorithm is said to be \texttt{convergent over }$\Omega$
iff for every $x^{0}\in\Omega$ there exist the limit $\lim_{k\rightarrow\infty}x^{k}=y\left(x^{0}\right)$
and $y\left(x^{0}\right)\in\Omega$. It is said to be \texttt{boundedly
convergent} \texttt{over }$\Omega$ iff, in addition, there exists
a $\gamma\geq0$ such that ${Prox}_{C}\left(y\left(x^{0}\right)\right)\leq\gamma$
for every $x^{0}\in\Omega$. \\

Next theorem, which gives sufficient conditions for strong perturbation
resilience of the Basic Algorithm, has been proved in \cite[Theorem 1]{hgdc12}
(in different wording).\\
\\
\textbf{Theorem 4.1}\label{theorem4.5} \emph{Assume that we are given
a $C\subseteq\mathbb{R}^{J}$, a proximity function ${Prox}_{C}$
and an algorithmic operator $\boldsymbol{A}_{C}:\mathbb{R}^{J}\rightarrow\Omega$.
If $\boldsymbol{A}_{C}$ is nonexpansive and is such that it defines
a boundedly convergent Basic Algorithm and if the proximity function
${Prox}_{C}$ is uniformly continuous, then the Basic Algorithm defined
by $\boldsymbol{A}_{C}$ is strongly perturbation resilient. }\\

Along with the $C\subseteq\mathbb{R}^{J}$, we look at the objective
function $\phi:\mathbb{R}^{J}\rightarrow\mathbb{R}$, with the convention
that a point in $\mathbb{R}^{J}$ for which the value of $\phi$ is
smaller is considered \textit{superior} to a point in $\mathbb{R}^{J}$
for which the value of $\phi$ is larger. The essential idea of the
SM is to make use of the perturbations of (\ref{eq:perturb}) to transform
a strongly perturbation resilient algorithm that seeks a constraints-compatible
solution for $C$ into one whose outputs are equally good from the
point of view of constraints-compatibility, but are superior (not
necessarily optimal) according to the objective function $\phi$.

This is done by producing from the Basic Algorithm another algorithm,
called its \textit{superiorized} version, that makes sure not only
that the $\beta_{k}v^{k}$ are bounded perturbations, but also that
$\phi\left(y^{k}+\beta_{k}v^{k}\right)\leq\phi\left(y^{k}\right)$,
for all $k$. To do so, we use the next concept, closely related to
the concept of ``descent direction\textquotedblright{}.\\
\\
\textbf{Definition 4.4} Given a function $\phi:\mathbb{R}^{J}\rightarrow\mathbb{R}$
and a point $y\in\mathbb{R}^{J}$, we say that a vector $d\in\mathbb{R}^{J}$
is \texttt{nonascending} \texttt{for }$\phi$\texttt{ at }$y$ iff
$\left\Vert d\right\Vert \leq1$ and there is a $\delta>0$ such that
\begin{equation}
\text{for all }\lambda\in\left[0,\delta\right]\text{ we have }\phi\left(y+\lambda d\right)\leq\phi\left(y\right).\label{eq:nonascend}
\end{equation}
\\

Obviously, the zero vector is always such a vector, but for superiorization
to work we need a sharp inequality to occur in (\ref{eq:nonascend})
frequently enough.

The Superiorized Version of the Basic Algorithm assumes that we have
available a summable sequence $\left\{ \eta_{\ell}\right\} _{\ell=0}^{\infty}$
of positive real numbers (for example, $\eta_{\ell}=a^{\ell}$, where
$0<a<1$) and it generates, simultaneously with the sequence $\left\{ y^{k}\right\} _{k=0}^{\infty}$
in $\Omega$, sequences $\left\{ v^{k}\right\} _{k=0}^{\infty}$ and
$\left\{ \beta_{k}\right\} _{k=0}^{\infty}$. The latter is generated
as a subsequence of $\left\{ \eta_{\ell}\right\} _{\ell=0}^{\infty}$,
resulting in a nonnegative summable sequence $\left\{ \beta_{k}\right\} _{k=0}^{\infty}$.
The algorithm further depends on a specified initial point $y^{0}\in\Omega$
and on a positive integer $N$. It makes use of a logical variable
called \textit{loop}\emph{. }The superiorized algorithm is presented
next by its pseudo-code.\medskip{}

\label{alg_super}\textbf{Superiorized Version of the Basic Algorithm}
\begin{enumerate}
\item \textbf{set} $k=0$
\item \textbf{set} $y^{k}=y^{0}$
\item \textbf{set} $\ell=-1$
\item \textbf{repeat}
\item $\qquad$\textbf{set} $n=0$
\item $\qquad$\textbf{set} $y^{k,n}=y^{k}$
\item $\qquad$\textbf{while }$n$\textbf{$<$}$N$
\item $\qquad$\textbf{$\qquad$set }$v^{k,n}$\textbf{ }to be a nonascending
vector for $\phi$ at $y^{k,n}$
\item $\qquad$\textbf{$\qquad$set} \emph{loop=true}
\item $\qquad$\textbf{$\qquad$while}\emph{ loop}
\item $\qquad\qquad\qquad$\textbf{set $\ell=\ell+1$}
\item $\qquad\qquad\qquad$\textbf{set} $\beta_{k,n}=\eta_{\ell}$
\item $\qquad\qquad\qquad$\textbf{set} $z=y^{k,n}+\beta_{k,n}v^{k,n}$
\item $\qquad\qquad\qquad$\textbf{if }$\phi\left(z\right)$\textbf{$\leq$}$\phi\left(y^{k}\right)$\textbf{
then }
\item $\qquad\qquad\qquad\qquad$\textbf{set }$n$\textbf{$=$}$n+1$
\item $\qquad\qquad\qquad\qquad$\textbf{set }$y^{k,n}$\textbf{$=$}$z$
\item $\qquad\qquad\qquad\qquad$\textbf{set }\emph{loop = false}
\item $\qquad$\textbf{set }$y^{k+1}$\textbf{$=$}$\boldsymbol{A}_{C}\left(y^{k,N}\right)$
\item $\qquad$\textbf{set }$k=k+1$ \medskip{}

\end{enumerate}
\textbf{Theorem 4.2}\label{theorem4.5-1} \emph{Any sequence $\left\{ y^{k}\right\} _{k=0}^{\infty}$,
generated by the Superiorized Version of the Basic Algorithm, satisfies
(\ref{eq:perturb}). Further, if, for a given $\varepsilon>0,$ the
$\varepsilon$-output $O\left(C,\varepsilon,\left\{ x^{k}\right\} _{k=0}^{\infty}\right)$
of the Basic Algorithm exists for every $x^{0}\in\Omega$, then every
sequence $\left\{ y^{k}\right\} _{k=0}^{\infty}$, generated by the
Superiorized Version of the Basic Algorithm, has an $\varepsilon^{\prime}$-output
$O\left(C,\varepsilon',\left\{ y^{k}\right\} _{k=0}^{\infty}\right)$
for every $\varepsilon^{\prime}>\varepsilon$.}\\

$\qquad$This theorem follows from the analysis of the behavior of
the Superiorized Version of the Basic Algorithm in \cite{hgdc12}.
In other words, the Superiorized Version produces outputs that are
essentially as constraints-compatible as those produced by the original
not superiorized algorithm. However, due to the repeated steering
of the process by lines 7 to 17 toward reducing the value of the objective
function $\phi$, we can expect that the output of the Superiorized
Version will be superior (from the point of view of $\phi$) to the
output of the original algorithm.

\section{A Computational Demonstration\label{sect:experimental}}

\subsection{The x-ray CT problem}

The fully-discretized model in the series expansion approach to the
image reconstruction problem of x-ray computerized tomography (CT)
is formulated in the following manner. A Cartesian grid of square
picture elements, called \textit{pixels}, is introduced into the region
of interest so that it covers the whole picture that has to be reconstructed.
The pixels are numbered in some agreed manner, say from 1 (top left
corner pixel) to $J$ (bottom right corner pixel).

The x-ray attenuation function is assumed to take a constant value
$x_{j}$ throughout the $j$th pixel, for $j=1,2,...,J$. Sources
and detectors are assumed to be points and the rays between them are
assumed to be lines. Further, assume that the length of intersection
of the $i$th ray with the $j$th pixel, denoted by $a_{j}^{i}$,
for $i=1,2,...,I,\;\; j=1,2,...,J$, represents the weight of the
contribution of the $j$th pixel to the total attenuation along the
$i$th ray.

The physical measurement of the total attenuation along the $i$th
ray, denoted by $b_{i}$, represents the line integral of the unknown
attenuation function along the path of the ray. Therefore, in this
fully-discretized model, the line integral turns out to be a finite
sum and the model is described by a system of linear equations
\begin{equation}
\sum_{j=1}^{J}x_{j}a_{j}^{i}=b_{i},\text{ \ }i=1,2,\ldots,I.\label{eqn:CT-syseqn}
\end{equation}
In matrix notation we rewrite (\ref{eqn:CT-syseqn}) as 
\begin{equation}
Ax=b,\label{eqn:CT-yAx}
\end{equation}
where $b\in\mathbb{R}^{I}$ is the \textit{measurement vector}, $x\in\mathbb{R}^{J}$
is the \textit{image vector}, and the $I\times J$ matrix $A=\left(a_{j}^{i}\right)$
is the \textit{projection matrix}. See \cite{HERM09}, especially
Section 6.3, for a complete treatment of this subject.

\subsection{The algorithms that we use}

In this section we describe the PSM and SM algorithms specifically
used in our demonstration. We applied both algorithms to solve the
fully-discretized model in the series expansion approach to the image
reconstruction problem of x-ray CT, formulated in the previous subsection
and represented by the optimization problem
\begin{equation}
\mathrm{minimize}\left\{ \phi(x)\mid Ax=b\text{ and }0\leq x\leq1\right\} .\label{eq: problem min phi x}
\end{equation}

The box constraints are natural for this problem: If $x_{j}$ represents
the linear attenuation coefficient, measured in cm\textsuperscript{-1},
at a medically-used x-ray energy spectrum in the $j$th pixel, then
the box constraints $0\leq x\leq1$ are reasonable for tissues in
the human body; see Table 4.1 of \cite{HERM09}. Hence, for the image
reconstruction problem of x-ray CT, we define $\Omega$ by
\begin{equation}
\Omega=\left\{ x\in\mathbb{R}^{J}\mid0\leq x\leq1\right\} .\label{eq:OMEGA}
\end{equation}
We note that this $\Omega$ is bounded.

The choice of $C$ in (\ref{prob:cons-min}) is of the type specified
in (\ref{eq:inersection}), with $L=I+1$, $C_{i}=\left\{ x\in\mathbb{R}^{J}\mid\left\langle a^{i},x\right\rangle =b_{i}\right\} $,
for $i=1,2,\ldots,I$ and $C_{I+1}=\Omega$. Furthermore, since in
the experiment reported below, we start with a specific image vector
$x\in\Omega$ and calculate from it the measurement vector $b\in\mathbb{R}^{I}$
using (\ref{eqn:CT-syseqn}), we know that $C$ is a nonempty subset
of $\Omega$, which is the requirement stated below (\ref{eq:inersection}).

For any such $C$, we define ${Prox_{C}}:\Omega\rightarrow\mathbb{R}_{+}$
by 
\begin{equation}
{Prox}_{C}(x)=\sqrt{\sum\limits _{i=1}^{I}\left(b_{i}-\left\langle a^{i},x\right\rangle \right)^{2}}.\label{eq:Prox}
\end{equation}
Note that this proximity function ${Prox}_{C}$ is uniformly continuous
and thus satisfies the condition stated for it in Theorem \ref{theorem4.5}.

Our choice for the objective function $\phi$ is the total variation
(TV) of the image vector $x.$ Denoting the $G\times H$ image array
$X$ ($GH=J$) obtained from the image vector $x$ by $X_{g,h}=x_{(g-1)H+h}$,
for $1\leq g\leq G$ and $1\leq h\leq H$, we use 
\begin{equation}
\phi\left(x\right)=\mathrm{TV}(X)=\sum\limits _{g=1}^{G-1}\sum\limits _{h=1}^{H-1}\sqrt{\left(X_{g+1,h}-X_{g,h}\right)^{2}+\left(X_{g,h+1}-X_{g,h}\right)^{2}}.\label{eq:TV}
\end{equation}

\subsubsection{The Projected Subgradient Method\label{sub:Projected-Subgradient-Method}}

We implemented the PSM with the choice of $C$ and the objective function
$\phi$ described above. We used the PSM recursion formula (\ref{eq:sgp})
and adopted a nonsummable diminishing step-length rule of the form
$t_{k}=\gamma_{k}/\left\Vert \phi^{\prime}\left(x^{k}\right)\right\Vert $,
where $\gamma_{k}\geq0,\,\,\,\lim_{k\rightarrow\infty}\gamma_{k}=0,\,\,\,\sum_{k=0}^{\infty}\gamma_{k}=\infty.$
\medskip{}

\textbf{The PSM Algorithm}\\
(P1) \textbf{Initialization}: Select a point $x^{0}\in\mathbb{R}^{J}$,
select integers $K$ and $M$, use two real number variables $\boldsymbol{curr}$
and $\boldsymbol{prev}$, and set $\boldsymbol{curr}=\phi\left(x^{0}\right)$
and $\boldsymbol{prev}=\,\boldsymbol{curr}$.\\
(P2) \textbf{Iterative step}: Given the current iterate $x^{k}$,
calculate the next one as follows:\\
(P2.1) Calculate a subgradient of $\phi$ at $x^{k},$ i.e., $\phi^{\prime}\left(x^{k}\right)\in\partial\phi\left(x^{k}\right)$,
a step-size $t_{k}=k^{-1/4}/\left\Vert \phi^{\prime}\left(x^{k}\right)\right\Vert {}_{2}$
and the vector
\begin{equation}
q^{k}=x^{k}-t_{k}\phi^{\prime}\left(x^{k}\right).
\end{equation}
\\
(P2.2) Calculate the next iterate as the projection of $q^{k}$ onto
$C$ by solving 
\begin{equation}
x^{k+1}=\arg\min_{x}\left\{ \frac{1}{2}\left\Vert x-q^{k}\right\Vert {}^{2}\mid Ax=b\text{ and }0\leq x\leq1\right\} .\label{eq:29-1}
\end{equation}
\\
(P2.3) If $\phi\left(x^{k+1}\right)\leq$ $\boldsymbol{curr}$, then
$\boldsymbol{curr}=\phi\left(x^{k+1}\right)$.\\
(P3) \textbf{Stopping rule}: If $k\,{mod}\, K=0$ (i.e., $k$ is divisible
by $K$), then:\\
If $\boldsymbol{prev}\,-\,\boldsymbol{curr}\,<\,\boldsymbol{prev}\,/\, M$
then stop. Otherwise, $\boldsymbol{prev}=\,\boldsymbol{curr}\,$ and
go to (P2).\medskip{}

That the PSM algorithm converges to a solution of (\ref{prob:cons-min})
follows from \cite[Subsection 3.2.3]{nestrov-book}, in particular,
from Theorem 3.2.2 therein, provided that $\phi$ is convex and locally
Lipschitz continuous and $C$ is closed and convex. The latter is
indeed the case for the $C$ in (\ref{eq: problem min phi x}). The
convexity of the $\phi$ of (\ref{eq:TV}) follows from the end of
the Proof of Proposition 1 in \cite{cp04}. Its Lipschitz continuity
on the whole space $\mathbb{R}^{J}$ follows from the fact that the
TV function can be rewritten as 
\begin{equation}
TV(X)=\sum_{g=1}^{G-1}\sum_{h=1}^{H-1}\left\Vert A_{g,h}X\right\Vert _{2}.\label{tv_new-1-1}
\end{equation}
where $A_{g,h}$ is a square matrix having only two nonzero rows,
with the first nonzero row containing only two nonzero elements 1
and $-1$ that correspond to the variables $X_{g+1,h}$ and $X_{g,h}$,
respectively, and the second nonzero row containing only two nonzero
elements 1 and $-1$ that correspond to the variables $X_{g,h+1}$
$X_{g,h}$, respectively.

In our implementation we solved problem (\ref{eq:29-1}), in step
(P2.2) above, by considering its dual
\begin{equation}
\mathrm{maximize}\left\{ f(\lambda)\mid\lambda\in\mathbb{R}^{I}\right\} ,\label{dualproblem-1}
\end{equation}
where 
\begin{equation}
\begin{array}{rcl}
f(\lambda) & = & \frac{1}{2}\left\Vert q^{k}-A^{T}\lambda-P_{C_{I+1}}\left(q^{k}-A^{T}\lambda\right)\right\Vert ^{2}-\frac{1}{2}\left\Vert q^{k}-A^{T}\lambda\right\Vert ^{2}\\
 & - & \left\langle \lambda,b\right\rangle +\frac{1}{2}\left\Vert q^{k}\right\Vert ^{2}.
\end{array}\label{eq:E3-1-1}
\end{equation}
The optimal point $x^{\ast k}$ of (\ref{eq:29-1}) is then 
\begin{equation}
x^{\ast k}=P_{C_{I+1}}\left(q^{k}-A^{T}\lambda^{\ast k}\right),
\end{equation}
where $\lambda^{\ast k}$ is the optimal solution of (\ref{dualproblem-1}).
To find $\lambda^{\ast k}$ we minimized $-f(\lambda)$ using the
Optimal Method of Nesterov \cite{Nesterov83}, as generalized by G\"{u}ler
\cite[p. 188]{gul93}, whose generic description for unconstrained
minimization of a convex function $\theta(\lambda)$, which is continuously
differentiable with Lipschitz continuous gradient, is as follows.

(N1) \textbf{Initialization:} Select a $\mu^{0}\in\mathbb{R}^{J}$,
a positive $\alpha_{-1}$ and put $\lambda^{-1}=\mu^{0}$, $\beta_{0}=1$
and $k=0$.\\
(N2) \textbf{Iterative Step:} Given $\lambda^{k-1}$, $\mu^{k}$,
$\alpha_{k-1}$ and $\beta_{k}$: \\
(N2.1) Calculate the smallest index $s\geq0$ for which the following
inequality holds 
\begin{equation}
\theta\left(\mu^{k}\right)-\theta\left(\mu^{k}-2^{-s}\alpha_{k-1}\nabla\theta\left(\mu^{k}\right)\right)\geq2^{-s-1}\alpha_{k-1}\left\Vert \nabla\theta\left(\mu^{k}\right)\right\Vert ^{2}.
\end{equation}
 \\
(N2.2) Calculate the next iterate by 
\begin{equation}
\alpha_{k}=2^{-s}\alpha_{k-1}\text{ and }\,\lambda^{k}=\mu^{k}-\alpha_{k}\nabla\theta\left(\mu^{k}\right),
\end{equation}
and update 
\begin{equation}
\beta_{k+1}=\left(\frac{1}{2}+\frac{1}{2}\sqrt{4\beta_{k}^{2}+1}\right),
\end{equation}
and 
\begin{equation}
\mu^{k+1}=\lambda^{k}+\frac{\beta_{k}-1}{\beta_{k+1}}\left(\lambda^{k}-\lambda^{k-1}\right).
\end{equation}
\medskip{}
When a stopping rule applies, then the point $\lambda^{k}$ is the
output of the method.

In the reported experiments, we used the starting points $x^{0}$
in the PSM Algorithm and $\lambda^{-1}=\mu^{0}$ in (N1) above to
be zero vectors. In the initialization step of the PSM Algorithm,
we selected $K=10$ and $M=5000$. In (N1), we chose $\alpha_{-1}=10$.

\subsubsection{The Superiorization Method\label{sub:Superiorization-Method-(SM)}}

Our selected choice for the operator $\boldsymbol{A}_{C}$ in the
Basic Algorithm as well as in the Superiorized Version of the Basic
Algorithm, as described in Section \ref{sect:SM}, is based on an
algebraic reconstruction technique (ART), see \cite[Chapter 11]{HERM09}.
Specifically, for $i=1,2,\ldots,I$, we define the operators $U_{i}:\mathbb{R}^{J}\rightarrow\mathbb{R}^{J}$
by
\begin{equation}
U_{i}(x)=x+\frac{b_{i}-\left\langle a^{i},x\right\rangle }{\left\Vert a^{i}\right\Vert ^{2}}a^{i}.\label{block operator}
\end{equation}
Defining the projection operator onto the unit box $\Omega$ by $Q:\mathbb{R}^{J}\rightarrow\Omega$ 

\begin{equation}
\left(Q(x)\right)_{j}=\left\{ \begin{array}{cc}
x_{j}, & \text{if}\;0\leq x_{j}\leq1,\\
0, & \text{if}\;\; x_{j}<0,\\
1, & \text{if}\;\;1<x_{j},
\end{array}\right.\label{eqn:CT-basis-1}
\end{equation}
for $j=1,2,...,J$, we specify the algorithmic operator $\boldsymbol{A}_{C}:\Omega\rightarrow\Omega$
by 
\begin{equation}
\boldsymbol{A}_{C}\left(x\right)=QU_{I}\cdots U_{2}U_{1}(x).\label{eq:perturb-1}
\end{equation}
Since the individual $U_{i}$s as well as the $Q$ are clearly nonexpansive
operators, the same is true for $\boldsymbol{A}_{C}$. 

By well-known properties of ART (see, for example, Sections 11.2 and
15.8 of \cite{HERM09}), the Basic Algorithm with this algorithmic
operator is convergent over $\Omega$ and, in fact, for every $x^{0}\in\Omega$,
the limit $y\left(x^{0}\right)$ is in $C$. It follows that, for
every $x^{0}\in\Omega$, ${Prox}_{C}\left(y\left(x^{0}\right)\right)=0,$
and so the Basic Algorithm is boundedly convergent. According to Theorem
\ref{theorem4.5}, this combined with the facts that $\boldsymbol{A}_{C}$
is nonexpansive and the proximity function ${Prox}_{C}$ is uniformly
continuous, implies that the Basic Algorithm defined by $\boldsymbol{A}_{C}$
is strongly perturbation resilient.

The following uses the convergence of the Basic Algorithm to an element
of $C$ and Theorem 2. Since for all $\varepsilon>0,$ the $\varepsilon$-output
$O\left(C,\varepsilon,\left\{ x^{k}\right\} _{k=0}^{\infty}\right)$
of the Basic Algorithm is defined for every $x^{0}\in\Omega$, we
also have that every sequence $\left\{ y^{k}\right\} _{k=0}^{\infty}$
generated by the Superiorized Version of the Basic Algorithm has an
$\varepsilon^{\prime}$-output $O\left(C,\varepsilon',\left\{ y^{k}\right\} _{k=0}^{\infty}\right)$
for every $\varepsilon^{\prime}>0$. This means that for the specific
type of $C$ that is used in our comparative study, the Superiorized
Version of the Basic Algorithm is guaranteed to produce an $\varepsilon'$-compatible
output for any $\varepsilon'>0$ and any initial point $y^{0}\in\Omega$.

The specific choices made when running the Superiorized Version of
the Basic Algorithm for our comparative study were the following.
We selected $\eta_{\ell}=0.999^{\ell}$, $y^{0}$ to be the zero vector
and $N=9$. All these choices we made are based on auxiliary experiments
(not included in this paper) that helped determine optimal parameters
for the data-set discussed in Subsection \ref{sub:Computational-results}.
In addition, we need to specify how the nonascending vector $v^{k,n}$
is selected in line 8 of the Superiorized Version of the Basic Algorithm.
We use the method specified in \cite{hgdc12} (especially Section
II.D, the paragraph following equation (12) and Theorem 2 in the Appendix).
Specifically, we define another vector $w$ and set $v^{k,n}$ to
be the zero vector if $\left\Vert w\right\Vert =0$ and $-\frac{w}{\left\Vert w\right\Vert }$
otherwise. The components of $w$ are computed by $w_{j}=\frac{\partial\phi}{\partial x_{j}}(y^{k,n})$
if the partial derivative can be calculated without a numerical difficulty
and $w_{j}=0$ otherwise, for $1\leq j\leq J$. Looking at (\ref{eq:TV})
we see that formally the partial derivative $w_{j}=\frac{\partial\phi}{\partial x_{j}}(y^{k,n})$
is the sum of at most three fractions; the phrase ``numerical difficulty''
in the previous sentence refers to the situation when in one of these
fractions the denominator has an absolute value less than $10^{-20}$.

\subsection{The computational result\label{sub:Computational-results}}

The computational work reported here was done on a single machine
using a single CPU, an Intel i5-3570K 3.4 Ghz with 16 GB RAM using
the SNARK09 software package \cite{SNARK09,kdh13}; the phantom, the
data, the reconstructions and displays were all generated within this
same framework. In particular, this implies that differences in the
reported reconstruction times are not due to the different algorithms
being implemented in different environments.

\begin{figure}
\centering\includegraphics[scale=0.5]{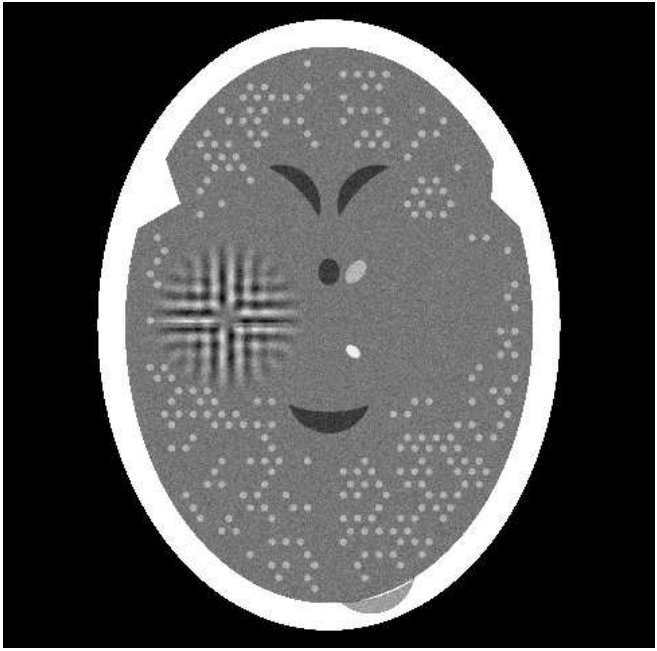}\caption{\label{fig:The-head-phantom.}The head phantom. The value of its TV
is 984. Its tomographic data was obtained for 60 views.}
\end{figure}

Figure \ref{fig:The-head-phantom.} shows the phantom used in our
study, which is a $485\times485$ digitized image whose TV is 984.
The phantom corresponds to a cross-section of a human head (based
on \cite[Figure 4.6]{HERM09}). It is represented by a vector with
$235,225$ components, each standing for the average x-ray attenuation
coefficient within a pixel. Each pixel is of size $0.376\times0.376$
mm$^{2}$. The values of the components are in the range of $[0,$
$0.6241749]$, however, the display range used here was much smaller,
namely $[0.204,$ $0.21675]$. The mapping between the two ranges
is such that any value below $0.204$ is shown as black and any value
above $0.21675$ is shown as white with a linear mapping in-between.
We used this display window for all images presented here.

\begin{figure}
\centering\subfloat[]{\includegraphics[scale=0.5]{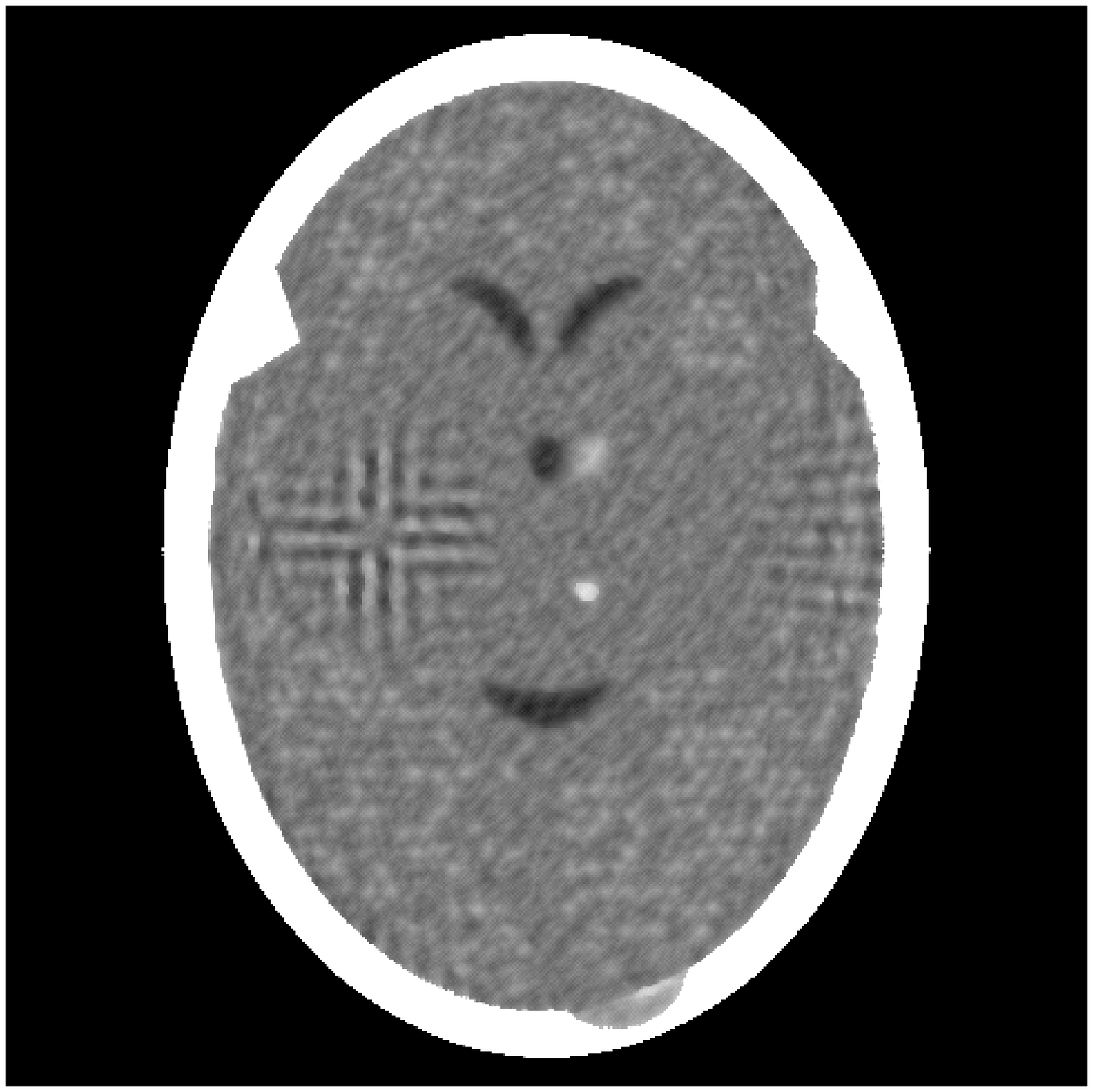}}\\
\subfloat[]{$\,$\includegraphics[scale=0.5]{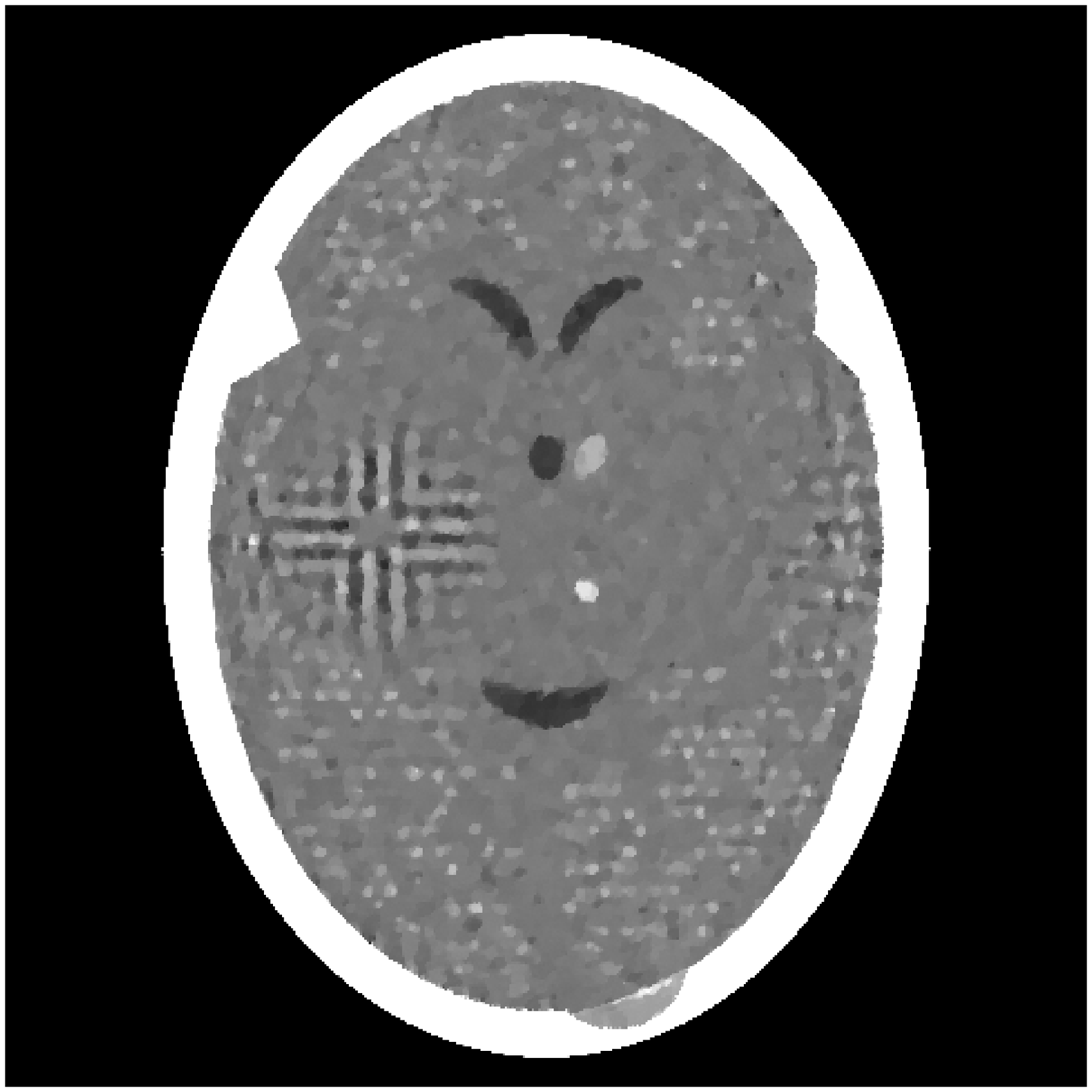}}

\caption{\label{fig:Reconatructions-of-the}Reconstructions of the head phantom
of Figure \ref{fig:The-head-phantom.}. (a) The image reconstructed
by the PSM has $TV=919$ and was obtained after 2217 seconds. (b)
The image reconstructed by the SM has $TV=873$ and was obtained after
102 seconds.}
\end{figure}

Data were collected by calculating line integrals through the digitized
head phantom in Figure \ref{fig:The-head-phantom.} using $60$ sets
of equally rotated (in $3$ degrees increments) parallel lines, with
lines in each set spaced at $0.752$ mm from each other. Each line
integral gives rise to a linear equation and represents a hyperplane
in $\mathbb{R}^{J}$. The phantom itself lies in the intersection
of all the hyperplanes that are associated with these lines, and it
also satisfies the box constraints in (\ref{eq:OMEGA}). The total
number of linear equations is $18,524$, making our problem underdetermined
with $235,225$ unknowns (the intersection of all the hyperplanes
is in an at least $216,701$-dimensional subspace of $R^{235,225}$).
In the comparative study, we first applied the PSM and then the SM
to these data as follows. 

The PSM was implemented as described in Subsection \ref{sub:Projected-Subgradient-Method}.
In particular, it started with the zero vector, for which ${Prox}_{C}\left(x^{0}\right)=326$.
It was stopped according to the Stopping Rule (P3), the iteration
number at that time was 815 and the value of the proximity function
was ${Prox}_{C}\left(x^{815}\right)=0.0422$, which is very much smaller
than the value at the initial point. The computer time required was
2217 seconds. The TV of the output was 919, which is less than that
of the phantom, indicating that the PSM is performing its task of
producing a constraints-compatible output with a low TV. This output
is shown in Figure \ref{fig:Reconatructions-of-the}(a).

We used the Superiorized Version of the Basic Algorithm, as described
in Subsection \ref{sub:Superiorization-Method-(SM)} to generate a
sequence $\left\{ y^{k}\right\} _{k=0}^{\infty}$ until it reached
$O\left(C,0.0422,\left\{ y^{k}\right\} _{k=0}^{\infty}\right)$ and
considered that to be the output of the SM. We know that this output
must exist for our problem and that its constraints-compatibility
will not be greater than that of the output of the PSM. The computer
time required to obtain this output was 102 seconds, which is over
twenty times shorter than what was needed by the PSM to get its output.
The TV of the the SM output was 876, which is also less than that
of the output of PSM. The SM output is shown in Figure \ref{fig:Reconatructions-of-the}(b).

\begin{table}
\centering%
\begin{tabular}{|c|c|c|}
\hline 
 & $TV$ value & Time (seconds)\tabularnewline
\hline 
\hline 
PSM & 919 & 2217\tabularnewline
\hline 
SM & 873 & $\;$102\tabularnewline
\hline 
\end{tabular}\caption{\label{tab:Performance-comparison-of}Performance comparison of the
PSM and the SM when producing the reconstructions in Figure \ref{fig:Reconatructions-of-the}.}
\end{table}

As summarized in Table \ref{tab:Performance-comparison-of}, with
the stopping rule that guarantees that the output of the SM is at
least as constraints-compatible as the output of the PSM, the SM showed
superior efficacy compared to the PSM: it obtained a result with a
lower TV value at less than one twentieth of the computational cost.

\section{Conclusions\label{sect:conclusions}}

The superiorization methodology (SM) allows the conversion of a feasibility-seeking
algorithm, designed to find an $\varepsilon$-compatible solution
of the constraints, into a superiorized algorithm that inserts, into
the feasibility-seeking algorithm, objective function reduction steps
while preserving the guaranteed feasibility-seeking nature of the
algorithm. The superiorized algorithm interlaces objective function
nonascent steps into the original process in an automatic manner.
In case of strong perturbation resilience of the original feasibility-seeking
algorithm, mathematical results indicate why the superiorized algorithm
will be efficacious for producing an $\varepsilon$-compatible solution
output with a low value of the objective function.

We have presented an example for which the SM finds a better solution
to a constrained minimization problems than the\textit{ }projected
subgradient method (PSM), and in significantly less computation time.
This finding is understandable in view of the nature of how the methods
interlace feasibility-oriented activities with optimization activities.
While the PSM requires a projection onto the feasible region of the
constrained minimization problem, the SM needs to do only projections
onto the individual constraints whose intersection is the feasible
region. We demonstrated this experimentally on a large-sized application
that is modeled, and needs to be solved, as a constrained minimization
problem.

\bigskip{}

\textbf{Acknowledgments}. We thank the editor and reviewer for their
constructive comments. We would like to acknowledge the generous support
by Dr. Ernesto Gomez and Dr. Keith Schubert in allowing us to use
the GPU cluster at the Department of Computer Science and Engineering
at California State University San Bernardino. We are also grateful
to Joanna Klukowska for her advice on using optimized compilation
for speeding up SNARK09. This work was supported by the United States-Israel
Binational Science Foundation (BSF) Grant No. 200912, the U.S. Department
of Defense Prostate Cancer Research Program Award No. W81XWH-12-1-0122,
the National Science Foundation Award No. DMS-1114901, the U.S. Department
of Army Award No. W81XWH-10-1-0170, and by Grant No. R01EB013118 from
the National Institute of Biomedical Imaging and Bioengineering and
the National Science Foundation. The contents of this publication
is solely the responsibility of the authors and does not necessarily
represent the official views of the National Institute of Biomedical
Imaging and Bioengineering or the National Institutes of Health. \medskip{}


\begin{thebibliography}{10}
\bibitem{lz12}Luo, S., Zhou, T.: Superiorization of EM algorithm
and its application in single-photon emission computed tomography
(SPECT). Inverse Problems and Imaging, accepted for publication
\end{thebibliography}
\end{document}